\documentclass{amsart}
\usepackage{latexsym}
\usepackage{amsmath}
\usepackage{amsthm}
\usepackage{amssymb}
\usepackage{hyperref}

\newtheorem{thm}{Theorem}[section]

\newtheorem{cor}[thm]{Corollary}
\newtheorem{defi}[thm]{Definition}

\newtheorem{prop}[thm]{Proposition}
\newtheorem{rk}[thm]{Remark}

\newenvironment{preuve}{\noindent {\it Proof}}{\hfill$\square$}

\def\R{\mathbb{R}}
\let \vphi=\varphi
\let\la=\lambda
\let \al=\alpha
\let \ga=\gamma
\let \s=\sigma
\def\eps{\varepsilon}
\newcommand{\vip}{\vskip.2cm}
\newcommand{\bea}{\begin{eqnarray}}
\newcommand{\eea}{\end{eqnarray}}
\newcommand{\beas}{\begin{eqnarray*}}
\newcommand{\eeas}{\end{eqnarray*}}
\newcommand{\PR}{\mathbb{P}}
\newcommand{\E}{\mathbb{E}}
\newcommand{\be}{\begin{equation}}
\newcommand{\ee}{\end{equation}}
\newcommand{\calc}{{\mathcal C}}
\newcommand{\cF}{{\mathcal F}}
\newcommand{\calp}{{\mathcal P}}
\newcommand{\loi}{\mathcal{ L}}
\newcommand{\ds}{\displaystyle}
\newcommand{\intot}{\displaystyle \int _0^t }
\newcommand{\intd}{\ds \int_0^t \!\!\int_{\R^3}}
\newcommand{\intdd}{\ds \int_0^t \!\!\int_{\R^3\times\R^3}}

\newcommand{\tf}{{\tilde f}}
\newcommand{\tg}{{\tilde g}}
\newcommand{\tv}{{\tilde v}}
\newcommand{\tz}{{\tilde z}}

\newcommand{\tV}{{\tilde V}}

\newcommand{\cJ}{{{\mathcal J}}}
\newcommand{\cW}{{{\mathcal W}}}

\newcommand{\vs}{{v_*}}

\begin{document}

\title[Well-posedness of the homogeneous Landau equation]
{Well-posedness of the spatially homogeneous Landau equation 
for soft potentials}

\author{Nicolas Fournier$^1$, H\'el\`ene Gu\'erin$^2$}

\footnotetext[1]{LAMA UMR 8050,
Facult\'e de Sciences et Technologies,
Universit\'e Paris Est, 61, avenue du G\'en\'eral de Gaulle, 94010 Cr\'eteil 
Cedex, France, {\tt nicolas.fournier@univ-paris12.fr}}

\footnotetext[2]{IRMAR UMR 6625,
Univ. Rennes 1, Campus de Beaulieu, 35042 Rennes Cedex, France,
\texttt{helene.guerin@univ-rennes1.fr}}

\def\shortauthorname{Nicolas Fournier, H\'el\`ene Gu\'erin}

\def\abstractname{Abstract}

\begin{abstract}
We consider the spatially homogeneous Landau equation of kinetic theory, and 
provide a differential inequality for the Wasserstein distance with quadratic 
cost between two solutions. We deduce some well-posedness results. The main 
difficulty is that this equation presents a singularity 
for small relative velocities.
Our uniqueness result is the first one in the important case of
soft potentials. Furthermore, it is almost optimal for a class of 
moderately soft potentials, that is for a moderate singularity.
Indeed, in such a case, our result applies for initial conditions with
finite mass, energy, and entropy. For the other moderatley 
soft potentials, we assume additionnally some moment 
conditions on the initial data. For very soft potentials, we obtain
only a local (in time) well-posedness result, under some integrability
conditions.
Our proof is probabilistic,
and uses a stochastic version of the Landau equation, in the spirit of
Tanaka \cite{tanaka}.
\end{abstract}

\maketitle

\smallskip

\smallskip

\section{Introduction and results}

\subsection{The Landau equation}

We consider the spatially homogeneous Landau equation in dimension $3$ for 
soft potentials. This equation of kinetic physics,
also called Fokker-Planck-Landau equation, has been derived from the Boltzmann 
equation when the grazing collisions prevail in the gas. It describes the 
density $f_t(v)$ of particles having the velocity $v\in\R^3$ at time 
$t\geq 0$ :

\begin{equation}\label{eq:Landau}
\partial_t f_t(v) =\frac{1}{2}\sum_{i,j=1}^3%
\partial_i \left\{
\int_{\R^3}a_{ij}( v-v^*) \Big[ f_t( v^*) \partial_j f_t( v) -f_t( v) 
\partial_j f_t( v^*) \Big]dv^* \right\},
\end{equation}
where $\partial_t =\frac{\partial }{\partial t}$, $\partial_i =
\frac{\partial }{\partial v_i}$ and  $a(z)$ is a symmetric nonnegative 
matrix, depending on a parameter $\gamma$ 
(we will deal here with soft potentials, that is $\gamma\in(-3,0)$),
defined by
\be \label{def:a}
a_{ij}(z)=|z|^{\ga}(|z|^{2}\delta_{ij}-z_iz_j).
\ee 

The weak form of (\ref{eq:Landau}) writes, for any test 
function $\vphi:\R^3\mapsto \R$ 
\beas 
\frac{d}{dt}\int_{\R^3}\!\vphi(v)f_t(v)dv=
\iint_{\R^3\times \R^3} f_t(dv)f_t(dv^*)L\vphi(v,v^*)
\eeas
where the operator $L$ is defined by 
\bea \label{def:L}
L\vphi(v,v^*)&=&\frac{1}{2}\sum_{i,j=1}^3a_{ij}(v-v^*)\partial_{ij}^2\vphi(v)
+\sum_{i=1}^3 b_i(v-v^*)\partial_i\vphi(v)\\
\label{def:b}
\text{with }\ b_i(z)&=&\sum_{j=1}^3\partial_j a_{ij}(z)
=-2|z|^\ga z_i,\text{ for }i=1,2,3.\eea
We observe that the solutions to (\ref{eq:Landau}) conserve, at least
formally, the mass, the momentum and the kinetic energy: for any $t\geq 0$,
$$\int_{\R^3} f_t(v)\vphi(v)dv=\int_{\R^3} f_0(v)\vphi(v)dv, \quad \text {for} 
\quad\vphi(v)=1,\, v,\, |v|^2.$$
We classically may assume without loss of generality that 
${\int_{\R^3} f_0(v)dv=1}$.

Another fundammental {\it a priori} estimate is the decay of entropy, 
that is solutions satisfy, at least formally, for all $t\geq 0$, 
$$\int_{\R^3} f_t(v)\log f_t(v) dv \leq \int_{\R^3} f_0(v)\log f_0(v) dv.$$

We refer to Villani \cite{villani:98b,villani:h} for many details
on this equation, its physical meaning, 
its derivation from the Boltzmann equation, and on what is known about
this equation.

\subsection{Existing results and goals}

One usually speaks of hard potentials for $\gamma>0$, Maxwell molecules
for $\gamma=0$, soft potentials for $\gamma\in (-3,0)$,
and Coulomb potential for $\gamma=-3$.

The Landau equation is a {\it continuous} approximation
of the Boltzmann equation: when there are infinitely many infinitesimally small
collisions, the particle velocities become continuous in time, which can
be modeled by equation (\ref{eq:Landau}). The most interesting case
is that of Coulomb potential, since then the Boltzmann equation seems
to be meaningless. 
Unfortunately, it is also the most difficult case to study.
However, the Landau equation can be derived from the Boltzmann equation
with {\it true very soft potentials}, that is $\gamma \in (-3,-1)$.
The main idea is that the more $\gamma$ is negative, the more the Landau 
equation is physically interesting. We refer again to \cite{villani:h}
for a detailed survey about such considerations.

Another possible issue concerns numerics for the Boltzmann equation
without cutoff: one can approximate grazing collisions by the Landau
equation.

\vip

Let us mention that existence of weak solutions, under physically
reasonnable assumptions on initial conditions, 
has been obtained by Villani \cite{villani:98b} for
all previously cited potentials. 

\vip

We study here the question of uniqueness (and stability with
respect to the initial condition). This question is of particular importance,
since uniqueness is needed to 
justify the derivation of the equation, the convergence
from the Boltzmann equation to the Landau equation, the convergence
of some numerical schemes, ...

Villani \cite{villani:98} has obtained uniqueness for
Maxwell molecules, and this was extended by Desvillettes-Villani 
\cite{desvillettes-villani:00} to the case of 
hard potentials. 
To our knowledge, there is no available result in the important case of 
soft potentials. We adapt in this paper the ideas of some recent
works on the Boltzmann equation, see 
\cite{fournier-mouhot,fournier-guerin} (see also Desvillettes-Mouhot 
\cite{desvillettes-mouhot} for other ideas). 
We will essentially prove here that uniqueness and stability hold in the 
following situations:

\vip

(a) for $\gamma \in (\gamma_0,0)$, with $\gamma_0=1-\sqrt 5 \simeq -1.236$, as
soon as $f_0$ satifies the sole physical assumptions of 
finite mass, energy and entropy;

(b) for $\gamma \in (-2,\gamma_0]$, as soon $f_0$ has finite mass, energy, 
entropy, plus some finite moment of order $q(\gamma)$ large enough;

(c) for $\gamma \in (-3,-2]$, as soon $f_0$ has a finite mass, energy,
and belongs to $L^p$, with $p>3/(3+\ga)$, and the result is local in time. 

\vip

Observe that (a) is extremely satisfying, and (c) is quite disappointing.

Observe also that we obtain some much better results than for the 
Boltzmann equation \cite{fournier-mouhot,fournier-guerin}, where well-posedness
is proved in the following cases: for $\ga \in (-0.61,0)$ and $f_0$ with finite
mass, energy, entropy; for $\ga \in (-1,-0.61)$ and $f_0$ with finite mass, 
energy, entropy, and a moment of order $q(\ga)$ sufficiently large, and for $\ga
\in (-3,-1]$ and $f_0 \in L^p$ with finite mass and energy, with 
$p>3/(3+\ga)$, and the result being local in time in the latest case.

\vip

In \cite{villani:98b}, Villani proved several results on the convergence
of the Boltzmann equation to the Landau equation. His results work
{\it up to extraction of a subsequence}. Of course, our uniqueness
result allows us to get a true convergence.

\subsection{Some notation}
Let us denote by $C^2_b$ the set of $C^2$-functions 
$\vphi:\R^3 \mapsto \R$ bounded with their derivatives of order $1$ and $2$. 
Let also $L^p(\R^3)$ be the space of measurable 
functions $f$ with
${ ||f||_{L^p}:=(\int_{\R^3}f^p(v)dv)^{1/p}<\infty}$, 
and let $\mathcal{P}(\R^3)$ be the space of probability measures on $\R^3$. 
For $k\geq 1$, we set
$$\mathcal{P}_k(\R^3)=\lbrace f\in\mathcal{P}(\R^3): m_k(f)<\infty \rbrace 
\quad \text{with}\quad m_k(f):=\int_{\R^3}|v|^kf(dv).$$
For $\al\in(-3,0)$, we introduce the space 
$\cJ_\alpha$ of probability measures $f$ on $\R^3$ such that
\begin{eqnarray}\label{cjalpha}
J_\alpha(f)&:=&\sup_{v\in\R^3} \int_{\R^3}|v-v_*|^\alpha f(dv_*) <\infty,
\end{eqnarray}
We denote by $L^\infty\big([0,T],\mathcal{P}_k\big)$,
$L^1\big([0,T],L^p\big)$, $L^1([0,T],\cJ_\alpha)$
the set of measurable families $(f_t)_{t\in[0,T]}$ of probability measures on 
$\R^3$ such that $\sup_{t\in[0,T]}m_k(f_t)<\infty$,  
$\int_0^T||f_t||_{L^p}dt<\infty$, $\int_0 ^T J_\alpha(f_t) dt < \infty$
respectively. Observe that (see \cite[(5.2)]{fournier-guerin}):

\begin{rk}\label{rk:Lp} 
For $\al\in (-3,0)$ and $p>3/(3+\al)$, there exists a constant 
$C_{\al,p}>0$ such that for all nonnegative measurable $f:\R^3 \mapsto\R$, 
any $v^*\in\R^d$ 
$$ J_\alpha(f) \leq ||f||_{L^1}+C_{\al,p}||f||_{L^p}.$$
\end{rk}

For a nonnegative function $f\in L^1(\R^3)$, we denote its entropy by
$$H(f)=\int_{\R^3}f(v)\log(f(v))dv.$$
Finally we denote $x\wedge y=\min(x,y)$,  $x\vee y=\max(x,y)$
for $x,y\in \R_+$, and $z.\tilde{z}$ the scalar product of 
$z,\tilde{z}\in\R^3$ and $\mathcal{L}(X)$ the distribution of a random variable $X$. For some set $A$ we write 
$\mathbf{1}_A$ the usual indicator function of $A$.

\subsection{Weak solutions}
We introduce now the notion of weak solution for the Landau equation.
Some much more refined definitions were introduced by Villani 
\cite{villani:98b} to allow solutions with only finite mass, energy,
and entropy.
\vip

We observe here that for $\vphi \in C^2_b$, $|L \vphi (v,v_*)|\leq C_\vphi
(|v-v_*|^{\gamma+1} + |v-v_*|^{\gamma+2} )$.
Thus if $\gamma \in [-1,0]$, $|L \vphi (v,v_*)|\leq C_\vphi (1+|v|^2+|v_*|^2)$,
while if $\gamma \in (-3,-1)$, 
$|L \vphi (v,v_*)|\leq C_\vphi (1+|v|^2+|v_*|^2+|v-v_*|^{\gamma+1})$.
This guarantees that all the terms are well-defined in the definition below.

\begin{defi}
Let $\ga \in (-3,0)$. Consider a measurable family 
$f=(f_t)_{t\in[0,T]}\in L^\infty\big([0,T],\mathcal{P}_2\big)$.
If $\ga \in (-3,-1)$, assume additionnaly that $f\in L^1([0,T],\cJ_{\ga+1})$.
We say that $f$ solves the Landau equation (\ref{eq:Landau}) if
for any $\vphi \in\mathcal{C}^2_b(\R^3)$ and any $t\in[0,T]$ 
\beas 
\int_{\R^3}\vphi(v)f_t(dv)= \int_{\R^3}\vphi(v)f_0(dv) + \intot ds
\iint_{\R^3\times\R^3} f_t(dv)f_t(dv^*)L\vphi(v,v^*),
\eeas
where the operator $L$ is defined by (\ref{def:L}).
\end{defi}

\subsection{The Wasserstein distance and the Monge-Kantorovich problem}

Given $\pi \in \mathcal{ P}_2(\R^{3}\times \R^3)$, we respectively denote by
$\pi_1$ and $\pi_2$ its first and second marginals on $\R^3$. For
two probability measures 
$\mu,\nu\in \mathcal{P}_2(\R^3)$ and 
$\pi \in \mathcal{ P}_2(\R^{3}\times \R^3)$, we write $\pi<^{\mu}_{\nu}$
if $\pi_1=\mu$ and $\pi_2=\nu$.

The Wassertein distance $\cW_2$ on $\mathcal{ P}_2(\R^3)$ is defined by
\beas
\cW_2^2(\mu,\nu)&=&\inf_{\pi<^{\mu}_{\nu}} \int_{\R^3\times \R^3} |x-y|^2 \pi(dx,dy)\\
&=&\inf\lbrace \E[|V-\tilde{V}|^2] : (V,\tilde{V})\in\R^{3}\times\R^3 
\text{ such that }\mathcal{L}(V)=\mu \text{ and } \mathcal{L}(\tilde{V})
=\nu\rbrace .
\eeas

The set $(\mathcal{ P}_2(\R^3),\cW_2)$ is a Polish space, see e.g. Rachev
and R\"uschendorf \cite{rachev:98}. The topology is stronger that
the usual weak topology (for more details, see Villani, \cite{villani:03} 
Theorem 7.12).

It is well-known (see e.g.  Villani \cite{villani:03} Chapter 1 for details)
that the infimum is actually a minimum, and that if $\mu$ (or $\nu$)
has a density with respect to the Lebesgue measure on $\R^3$, then there is
a unique $\pi<^{\mu}_{\nu}$ such that 
$\cW_2^2(\mu,\nu)=\int_{\R^3\times \R^3} |x-y|^2 \pi(dx,dy)$.
Moreover, if we consider a family $ (\mu_\la,\nu_\la)_{\la\in\Lambda}$ of 
$ \mathcal{ P}_2(\R^3)\times\mathcal{ P}_2(\R^3)$ such that the map 
$\la\mapsto (\mu_\la,\nu_\la)$ is measurable, and if $\mu_\lambda$
has a density for all $\lambda$, then the function 
$\la\mapsto \pi_\la$ is measurable, see Fontbona-Gu\'erin-M\'el\'eard 
\cite{fontbona:08}.


\subsection{A general inequality}
\setcounter{equation}{0}

Our main result reads as follows.

\begin{thm}\label{thm:uniqueness}
Let $T>0$ and $\ga\in(-3,0)$. Consider two weak solutions $f$ and 
$\tilde f$ to the Landau equation (\ref{eq:Landau}) such that 
$f,\tilde f\in L^\infty\big([0,T],\mathcal{P}_2\big)
\cap L^1([0,T],\mathcal{J}_\ga)$. We also assume that $f_t$ (or $\tf_t$)
has a density with respect the Lebesgue measure on $\R^3$ 
for each $t\in [0,T]$. 
Then there exists an explicit constant $C_\ga>0$, depending
only on $\gamma$, such that for all $t\in [0,T]$,
$$\cW^2_2(f_t,\tilde f_t)\leq \cW^2_2(f_0,\tilde f_0)
\exp\left(C_\ga\int_0^t(J_\ga(f_s)+J_\ga(\tilde f_s))dt\right).$$
\end{thm}

Observe that if $H(f_0)<\infty$, then $H(f_t)<\infty$ for all $t\geq 0$,
so that $f_t$ has a density for all times. Thus this result
always applies for solutions with finite entropy.

This result asserts that uniqueness and stability hold in 
$L^\infty\big([0,T],\mathcal{P}_2\big)\cap L^1([0,T],\mathcal{J}_\ga)$.
Using Remark \ref{rk:Lp}, we immediately deduce that  
uniqueness and stability also hold in 
$L^\infty\big([0,T],\mathcal{P}_2\big)\cap L^1([0,T],L^p)$,
as soon as $p>3/(3+\gamma)$.

\subsection{Applications}


We now show the existence of solutions in
$L^\infty\big([0,T],\mathcal{P}_2\big)\cap L^1([0,T],L^p)$.

\begin{cor}\label{thm:existence}
(i) Assume that $\ga\in(-2,0)$. Let $q(\gamma):=\ga^2/(2+\ga)$.
Let $f_0 \in \calp_2(\R^3)$ satisfy also $H(f_0)<\infty$ and
$m_q(f_0)<\infty$,
for some $q>q(\ga)$. Consider $p \in (3/(3+\ga),(3q-3\ga)/(q-3\ga))
\subset (3/(3+\ga),3)$.
Then the Landau equation (\ref{eq:Landau}) 
has an unique weak solution in 
$L^\infty\big([0,T],\mathcal{P}_2\big)\cap L^1([0,T],L^p)$.

(ii)  Assume that $\ga\in(-3,0)$, and let $p>3/(3+\ga)$. Let
$f_0 \in \calp_2(\R^3)\cap L^p(\R^3)$.  Then there exists $T^*>0$ 
depending on $\ga, p,||f_0||_p$ such that there is an unique 
weak solution in $L^\infty_{loc}([0,T_*),L^p)$.
\end{cor}

Observe that if $\gamma \in (1-\sqrt 5,0)\simeq (-1.236,0)$, then
$q(\ga)<2$, and thus we obtain the well-posedness for the Landau
equation under the physical assumptions of finite mass,
energy, and entropy. This is of course extremely satisfying.

\subsection{Plan of the paper}
The paper is organized as follows: we first prove Theorem \ref{thm:uniqueness} 
in Section
\ref{ineq}, by means of stochastic Landau processes.
Next, we prove Corollary \ref{thm:existence} in Section \ref{appli}.

\section{A general inequality}\label{ineq}

We give in this section the proof of Theorem \ref{thm:uniqueness}.
In the whole section, $T>0$ and $\gamma \in (-3,0)$ are fixed,
and we consider two weak solutions $f=(f_t)_{t\in [0,T]}$ and 
$\tilde f=(\tf_t)_{t\in [0,T]}$ to (\ref{eq:Landau}) belonging to 
$L^\infty\big([0,T],\mathcal{P}_2\big)\cap L^1([0,T],\mathcal{J}_\ga)$.
We also assume that $f_t$ has a density with respect to the Lebesgue measure 
on $\R^3$ for each $t\geq 0$, which implies the uniqueness of the minimizer 
in $\cW_2(f_t,\tf_t)$.

\vip

We now introduce two coupled Landau stochastic processes,
the first one associated with $f$, the second one associated with 
$\tf$, in such a way
that they remain as close to each other as possible. 
Our approach is inspired by  \cite{guerin:03}, which was itself
inspired by the work of Tanaka \cite{tanaka} on the Boltzmann equation.

\vip

For any $s\in[0,T]$, we denote by $R_s$ the unique solution of the 
Monge-Kantorovich transportation problem for the couple $(f_s,\tilde f_s)$. 
Recall that $R_s(dv,d\tv)$ is a probability measure on $\R^3\times \R^3$,
with marginals $f_s$ and $\tf_s$, and that 
$\cW_2^2(f_s,\tf_s)=\iint_{\R^3\times\R^3} |v-\tv|^2R_s(dv,d\tv)$.
Let us notice that the map $s \mapsto R_s$ is measurable thanks to 
Theorem 1.3 of \cite{fontbona:08}. 

\vip
On some probability space $(\Omega,\cF,\PR)$,
we consider $W(dv,d\tv,ds)=(W_1,W_2,W_3)(dv,d\tv,ds)$ a three-dimensional 
space-time white noise on 
$\R^3\times \R^3\times[0,T]$ with covariance measure 
$R_s(dv,d\tilde v)ds$ (in the sense of Walsh \cite{walsh:84}): 
$W$ is an orthogonal martingale measure with covariance $R_s(dv,d\tilde v)ds$. 

Let us now consider two random variables $V_0,\tilde V_0$ with valued in 
$\R^3$ with laws $f_0,\tf_0$, independant of $W$, such that 
$\cW_2^2(f_0,\tilde f_0)=\E[|V_0-\tilde V_0|^2]$.
We finally consider the two following $\R^3$-valued 
stochastic differential equations.
\bea\label{eq:eds}
V_t&=&V_0+\intdd \sigma(V_s-v)W(dv,d\tilde v,ds)+\intd b(V_s-v)f_s(dv)ds,\\
\tilde V_t&=&\tilde V_0+\intdd \sigma(\tilde V_s-\tilde v)W(dv,d\tilde v,ds)
+\intd b(\tilde V_s-\tilde v)f_s(d \tilde v)ds, \label{eq:edstilde}
\eea
where $\s$ is a square root matrix of $a$
(recall (\ref{def:a})), that is $a(z)=\s(z)\s^t(z)$,
defined by
\begin{equation}\label{def:s}
\sigma \left( z\right) =\left| z\right| ^{\frac{\gamma }{2}}\left[ 
\begin{array}{ccc}
z_{2} & -z_{3} & 0 \\ 
-z_{1} & 0 & z_{3} \\ 
0 & z_{1} & -z_{2}
\end{array}
\right].
\end{equation}
We denote $(\mathcal{F}_t)_{t\geq 0}$ the filtration generated by $W$ and 
$V_0,\tV_0$, 
that is $\mathcal{F}_t=\s\{W([0,s]\times A), \;  
s\in[0,t], A \in {\mathcal B}(\R^3\times\R^3)\}\vee \mathcal{N}$,
where $\mathcal{N}$ is  the set of all $\mathbb{P}$-null subsets.

\begin{prop}\label{prop:edslin} 
(i) There exists a unique pair of continuous  $(\mathcal{F}_t)_{t\geq 0}$-adapted
processes $(V_t)_{t\in [0,T]}$, $(\tV_t)_{t\in [0,T]}$ solutions to
(\ref{eq:eds}) and (\ref{eq:edstilde}), satisfying 
$\E[\sup_{[0,T]}(|V_t|^2+|\tV_t|^2)]<\infty$.

(ii) For $t\in [0,T]$, consider the distributions $g_t=\loi(V_t)$ and 
$\tg_t=\loi(\tV_t)$. Then $g_t=f_t$ and $\tg_t=\tf_t$ for all
$t\in [0,T]$.
\end{prop}

The processes $(V_t)_{t\in [0,T]}$ and $(\tV_t)_{t\in [0,T]}$ 
are called \emph{Landau processes} 
associated with $f$ and $\tilde f$ respectively. 
We start the proof with a simple remark.

\begin{rk}\label{rk:estimsb}
Let us consider the functions $\s$ and $b$ respectively defined by 
(\ref{def:s}) and (\ref{def:b}), and recall that $\ga<0$. 
There exists a constant $C_\ga>0$ such that for any $z,\tilde z\in \R^3$,
\beas
|\s(z)-\s(\tilde z)|\leq& C_\ga|z-\tilde z|\left(|z|^{\frac{\ga}{2}}+
|\tilde z|^{\frac{\ga}{2}}\right)
\quad\text{and} \quad|b(z)-b(\tilde z)|\leq  C_\ga|z-\tilde z|
\left(|z|^{\ga}+|\tilde z|^{\ga}\right).
\eeas
\end{rk}

\begin{proof}
We have
\beas
|\s(z)-\s(\tilde z)|&\leq&\left||z|^{\ga/2}-|\tilde z|^{\ga/2}\right|\,|z|
+|z-\tilde z|\,|\tilde z|^{\ga/2}\\
&\leq &\frac{|\ga|}{2} |z|\,|z-\tilde z|
\max(|z|^{\frac{\ga}{2}-1},|\tilde z|^{\frac{\ga}{2}-1})
+(|z|^{\frac{\ga}{2}}+|\tz|^{\frac{\ga}{2}})|z-\tilde z|.
\eeas
By symmetry, we deduce that
\beas
|\s(z)-\s(\tilde z)| 
&\leq& |z-\tilde z|\left(\frac{|\ga|}{2} 
\min(|z|,|\tz|) \max(|z|^{\frac{\ga}{2}-1},|\tilde z|^{\frac{\ga}{2}-1})
+|\tilde z|^{\frac{\ga}{2}}+|z|^{\frac{\ga}{2}}\right)\\
&\leq&\left( \frac{|\ga|}{2}+1 \right)|z-\tilde z|
\left(|z|^{\frac{\ga}{2}}+|\tilde z|^{\frac{\ga}{2}}\right)\,
\eeas
The same computation works with $b$, starting with
$|b(z)-b(\tz)| \leq 2\left||z|^{\ga}-|\tilde z|^{\ga}\right|\,|z|
+2|z-\tilde z|\,|\tilde z|^{\ga}$.
\end{proof}

\begin{preuve} {\it of Proposition \ref{prop:edslin}.}
We just deal with (\ref{eq:eds}), the study of (\ref{eq:edstilde}) 
being of course the same.

\vip

{\it Step 1.} We start with the proof of (i), that is 
existence and uniqueness of a solution for (\ref{eq:eds}). 
We consider the map $\Phi$ which associates to a continuous adapted 
process $X=(X_t)_{t\in [0,T]}\in\calc([0,T],\R^3)$, with 
$\E[\sup_{0\leq s\leq T}|X_s|^2]<\infty$, the continuous adapted process
$\Phi(X)\in\calc([0,T],\R^3)$ defined by
$$ 
\Phi(X)_t=V_0+\intdd \sigma(X_s-v)W(dv,d\tilde v,ds)+\intd b(X_s-v)f_s(dv)ds.
$$ 
{\it Step 1.1.} Let us first prove that
\begin{equation*}
\E\left[\sup_{[0,T]} |\Phi(X)_t|^2\right] \leq C_T\left(1+\E[|V_0|^2]
+ \E[\sup_{[0,T]} |X_t|^2]+\sup_{[0,T]}m_2(f_s)
+\left(\int_0^T \!\!J_\ga(f_s)ds\right)^2 \right),
\end{equation*}
which is finite thanks to the conditions imposed on $V_0$, $f$, and $X$.
Using Doob's inequality, we easily get, for some constant $C$,
\beas
\E\left[\sup_{[0,T]} |\Phi(X)_t|^2\right] &\leq &
C\E[|V_0|^2] + C \int_0^T \iint_{\R^3\times\R^3} 
\E[|\sigma(X_s-v)|^2]R_s(dv,d\tv)ds \\
&&+ C\E\left[\left(\int_0^T \int_{\R^3} | b(X_s-v)|f_s(dv)ds\right)^2\right]
=:C(\E[|V_0|^2]+A+B).
\eeas
Using that the first marginal of $R_s$ is $f_s$, that 
$|\s(z)|^2\leq |z|^{\ga+2}\leq (1+|z|^2+|z|^\gamma)$, we observe that
\beas
A&\leq&  C \int_0^T \int_{\R^3} (1+\E[|X_s|^2]+|v|^2+ \E[|X_s-v|^\ga] f_s(dv)ds\\
&\leq & C \int_0^T (1+\E[|X_s|^2] + m_2(f_s)+J_\ga(f_s))ds,
\eeas
by definition of $J_\ga$, see (\ref{cjalpha}). Next, using 
that $|b(z)|=2|z|^{\ga+1}\leq 2+2|z|^\gamma+2|z|$,
\beas
B &\leq &  C \E\left[\left( \int_0^T \int_{\R^3} (1+|X_s|+|v|+|X_s-v|^\ga) 
f_s(dv)ds \right)^2\right] \\
&\leq & C \left[ T^2+ T^2 \E[\sup_{[0,T]} |X_s|^2] + T \int_0^T m_2(f_s) ds 
+ \left( \int_0^T J_\ga (f_s) ds \right)^2 \right].
\eeas
{\it Step 1.2.} 
Let us now consider two adapted processes $X,Y \in\calc([0,T],\R^3)$,
and show that
\beas
\E\left[\sup_{[0,t]} |\Phi(X)_s-\Phi(Y)_s|^2\right] \leq 
C_\ga  \left(1+ \intot J_\ga(f_s)ds\right) \intot \E[|X_s-Y_s|^2] J_\ga(f_s)ds.
\eeas
Arguing as previously and using Remark \ref{rk:estimsb}, we deduce that
\begin{align*}
\E\big[\sup_{[0,t]} |\Phi(X)_s&-\Phi(Y)_s|^2\big] \leq  
2 \intot \int_{\R^3\times\R^3}\E \left[|\sigma(X_s-v)-\sigma(Y_s-v)|^2\right] 
R_s(dv,d\tv) ds \\
&+ 2\E\left[\left(\intot \int_{\R^3} |b(X_s-v)-b(Y_s-v)| f_s(dv)ds \right)^2
\right]\\
&\leq C_\ga \intot \int_{\R^3} \E\left[|X_s-Y_s|^2(|X_s-v|^\ga+|Y_s-v|^\ga ) 
\right] f_s(dv)ds \\
&+ C_\ga  \E\left[ \left(\intot \int_{\R^3} |X_s-Y_s|(|X_s-v|^\ga+|Y_s-v|^\ga )
f_s(dv)ds  \right)^2\right] \\
&\leq C_\ga \intot \E[|X_s-Y_s|^2] J_\ga(f_s)ds 
+ C_\ga \E\left[ \left( \intot |X_s-Y_s| J_\ga(f_s)ds \right)^2\right].
\end{align*}
We conclude by using the Cauchy-Schwarz inequality.

{\it Step 1.3.} The uniqueness of the solution to (\ref{eq:eds}) immediately
follows from Step 1.2. Indeed, consider two solutions $V^1$ and $V^2$.
Then $V^1=\Phi(V^1)$ and $V^2=\Phi(V^2)$, so that Step 1.2. implies
that $\E[\sup_{[0,t]}|V^1_s-V^2_s|^2]
\leq C\int_0^t E[|V^1_s-V^2_s|^2] J_\ga(f_s) ds$ for some constant $C>0$ 
depending on $\ga,f$. Since
$t\mapsto J_\ga(f_t) \in L^1([0,T])$ by assumption, 
the Gronwall Lemma implies that
$\E[\sup_{[0,t]}|V^1_t-V^2_t|^2]=0$, whence $V^1=V^2$.

{\it Step 1.4.} Finally, one classically obtains the existence of
a solution using a Picard iteraction: consider the process
$V^0$ defined by $V^0_t=V_0$, and then define by induction $V^{n+1}=\Phi(V^n)$
(this is well-defined thanks to Step 1.1.). 
Using Step 1.2 and classical arguments, one easily checks that there exists
a continuous adapted process $V$ such that $\E[\sup_{[0,T]}|V_t-V^n_t|^2]$
tends to $0$. It is not difficult to pass to the limit in $V^{n+1}=\Phi(V^n)$,
whence $V=\Phi(V)$, and thus $V$ solves (\ref{eq:eds}).

\vip

{\it Step 2.} We now prove (ii). Let $V$ be the unique solution of 
(\ref{eq:eds}) and $g_s=\loi(V_s)$ for all $s\in [0,T]$.

{\it Step 2.1.} We first check that the family $g$ solves the 
\emph{linear Landau equation}: for any $\vphi\in\calc^2_b(\R^3)$,
\begin{equation}\label{eq:lineq}
\int_{\R^3}\vphi(x)g_t(dx)= \int_{\R^3}\vphi(x)f_0(dx)+
\intot \iint_{\R^3\times \R^3} L\vphi(x,v) g_s(dx)f_s(dv) ds .
\end{equation}
with $L$ defined by (\ref{def:L}). Applying the It\^o formula, we immediately
get
\beas
\vphi(V_t)=\vphi(V_0)+ \intot \int_{\R^3\times \R^3} \sum_{i,j} \partial_i 
\vphi(V_s) \sigma_{i j}(V_s-v) W_j(dv,d\tv,ds) \\
+\intot \int_{\R^3} \sum_i \partial_i \vphi(V_s) b_i(V_s-v) f_s(dv)ds\\
+ \frac{1}{2} \intot \int_{\R^3\times \R^3} \sum_{i,j,k} \partial_{ij} 
\vphi(V_s) \sigma_{ik}(V_s-v) \sigma_{jk}(V_s-v) R_s(dv,d\tv)ds.\\
\eeas
Taking expectations (which makes vanish the first integral), using that
the first marginal of $R_s$ is $f_s$ and that $\loi(V_0)=f_0$, we obtain
\beas
\int_{\R^3} \vphi(x) g_t(dx)= \int_{\R^3} \vphi(x) f_0(dx)
+\intot \int_{\R^3\times \R^3} \sum_i \partial_i \vphi(x) b_i(x-v) g_s(dx)f_s(dv)
ds\\
+ \frac{1}{2} \intot \int_{\R^3\times \R^3} \sum_{i,j,k} \partial_{ij} 
\vphi(x) \sigma_{ik}(x-v) \sigma_{jk}(x-v) g_s(dx)f_s(dv)ds,
\eeas
from which the conclusion follows, recalling (\ref{def:L}), 
since $\sigma.\sigma^t=a$.

{\it Step 2.2.}
One may apply the general uniqueness result of Bhatt-Karandikar
\cite[Theorem 5.2]{bhatt:93}, which implies uniqueness for 
(\ref{eq:lineq}). We omit
the proof here, but we refer to \cite[Lemma 4.6]{fournier-guerin} 
for a very similar proof concerning the Boltzmann equation. 
The main idea is that roughly, Bhatt-Karandikar have proved that
existence and uniqueness for a stochastic differential equation implies
uniqueness for the associated Kolmogorov equation, so that essentially,
point (i) implies uniqueness for (\ref{eq:lineq}).

{\it Step 2.3.} But $f$, being a weak solution to (\ref{eq:Landau}),
is also a weak solution to (\ref{eq:lineq}). We deduce that for all
$t\in [0,T]$, $g_t=f_t$.
\end{preuve}

\vip

To prove Proposition \ref{prop:edslin} there was no need to couple the 
stochastic processes $V$ and $\tilde V$ with the same white noise. 
But to evaluate the Wasserstein distance between two Landau solutions $f$ 
and $\tilde f$ using the stochastic processes, it is essential to connect 
them with the same white noise as we can see below.

\begin{prop}\label{prop:eval}
Consider the unique solutions $V$ and $\tV$ to (\ref{eq:eds}) and 
(\ref{eq:edstilde}) defined in Proposition \ref{prop:edslin}.
There exists a constant $C_\ga>0$ depending only on $\ga$ such
that
$$
\E[|V_t-\tilde V_t|^2]\leq \E[|V_0-\tilde V_0|^2]+ C_\ga\intdd \! \! \!\!
\E\big[ \big(|V_s-\tilde V_s|^2+ |v-\tilde v|^2\big)
\big(|V_s-v|^{\ga}+|\tilde V_s-\tilde v|^{\ga} \big)  \big]R_s(dv,d\tilde v) ds.
$$
\end{prop}

\begin{proof} First of all, we observe that since $R_s(dv,d\tv)$ has the
marginals $f_s(dv)$ and $\tf_s(d\tv)$, we may rewrite equations 
(\ref{eq:eds}) and (\ref{eq:edstilde}) as
\beas
V_t&=&V_0+\intdd \sigma(V_s-v)W(dv,d\tilde v,ds)+\intdd b(V_s-v)R_s(dv,d\tv)ds,
\\
\tilde V_t&=&\tilde V_0+\intdd \sigma(\tilde V_s-\tilde v)W(dv,d\tilde v,ds)
+\intdd b(\tilde V_s-\tilde v) R_s(dv,d\tv)ds.
\eeas
Using the It\^o formula and taking expectations, we obtain
\beas
\E[|V_t-\tilde V_t|^2]&=&\E[|V_0-\tilde V_0|^2] +\sum_{i,l=1}^3 \intdd 
\E\big[ \big(\s_{il}(V_s-v)-\s_{il}(\tilde V_s-\tilde v)\big)^2\big]
R_s(dv,d\tilde v) ds \\
&& +2\intdd \E\big[ \big(b(V_s-v)-b(\tilde V_s-\tilde v)\big).
(V_{s}-\tilde V_{s})\big]R_s(dv,d\tilde v) ds .
\eeas
Using Remark \ref{rk:estimsb}, we deduce that for some constant $C_\ga$,
\beas
\E[|V_t-\tilde V_t|^2]
&\leq&\E[|V_0-\tilde V_0|^2] \\
& + &C_\ga\intdd \E\big[ |V_s-\tilde V_s-v+\tilde v|^2(|V_s-v|^{\ga}
+|\tilde V_s-\tilde v|^{\ga})  \big]R_s(dv,d\tilde v) ds \\
& +&C_\ga\intdd \E\big[|V_s-\tilde V_s-v+\tilde v||V_s-\tilde V_s|
(|V_s-v|^{\ga}+|\tilde V_s-\tilde v|^{\ga}) \big]R_s(dv,d\tilde v) ds, \\
\eeas
from which the result immediately follows.
\end{proof}

\vip

We are finally in a position to conclude this section.

\vip

\begin{preuve} {\it of Theorem \ref{thm:uniqueness}.}
Let us recall briefly the situation. We have two weak solutions $f$ and 
$\tilde f$ to the Landau equation, belonging to 
$L^\infty\big([0,T],\mathcal{P}_2\big)\cap L^1([0,T],\mathcal{J}_\ga)$.
For each $s\in [0,T]$, $R_s$ has marginals $f_s,\tf_s$, and
satisfies  $\cW_2^2(f_s,\tilde f_s)=\iint_{\R^3\times\R^3}|v-\tilde v|^2
R_s(dv,d\tilde v)$. Then we have introduced the solutions $V$ and
$\tV$ to (\ref{eq:eds}-\ref{eq:edstilde}), and we have shown that
for each $t\in [0,T]$, $\loi(V_t)=f_t$ and $\loi(\tV_t)=\tf_t$.
An immediate consequence of this is that 
$\cW^2_2(f_t,\tf_t)\leq \E[|V_t-\tV_t|^2]$.

Using Proposition \ref{prop:eval}, we get
\beas
\E[|V_t-\tilde V_t|^2]&\leq& \E[|V_0-\tilde V_0|^2]+ C_\ga\intdd 
\E[|V_s-\tilde V_s|^2]  \big(|V_s-v|^{\ga}+
|\tilde V_s-\tilde v|^{\ga} \big)  \big]R_s(dv,d\tilde v)ds \\
&& +C_\ga\intdd  |v-\tilde v|^2\E\big(|V_s-v|^{\ga}+
|\tilde V_s-\tilde v|^{\ga} \big) R_s(dv,d\tilde v)ds.
\eeas
But since $\loi(V_s)=f_s$, we have $\E(|V_s-v|^{\ga})=\int_{\R^3}|x-v|^\ga f_s(dx)
\leq J_\ga(f_s)$. By the same way, $\E(|\tV_s-\tv|^{\ga})\leq J_\ga(\tf_s)$.
On the other hand, since the first marginal of $R_s$ is $f_s$, we deduce
that $\int_{\R^3\times\R^3} |V_s-v|^{\ga} R_s(dv,d\tv) =  
\int_{\R^3} |V_s-v|^{\ga} f_s(dv) \leq J_\ga(f_s)$, and by the same way,
 $\int_{\R^3\times\R^3} |\tV_s-\tv|^{\ga} R_s(dv,d\tv) \leq J_\ga(\tf_s)$.
Thus, since $\E[|V_0-\tilde V_0|^2]=\cW_2^2(f_0,\tilde f_0)$
and since $\cW^2_2(f_t,\tf_t)\leq \E[|V_t-\tV_t|^2]$,
\begin{align*}
\E[|V_t-\tilde V_t|^2]
&\leq\cW_2^2(f_0,\tilde f_0) + C_\ga\intot (J_\ga(f_s)+J_\ga(\tilde f_s))
\bigg[\E[|V_s-\tilde V_s|^2]+\int_{\R^3\times\R^3}\!\!\! |v-\tilde v|^2  
R_s(dv,d\tilde v)\bigg]ds\\
&\leq\cW_2^2(f_0,\tilde f_0)+C_\ga\intot (J_\ga(f_s)+J_\ga(\tilde f_s))
\bigg[\E[|V_s-\tilde V_s|^2]+\cW_2^2(f_s,\tilde f_s)\bigg]ds\\
&\leq \cW_2^2(f_0,\tilde f_0)+C_\ga\intot (J_\ga(f_s)+J_\ga(\tilde f_s))
\E[|V_s-\tilde V_s|^2]ds.
\end{align*}
Using finally the Gronwall Lemma, we get
$$ 
\E[|V_t-\tilde V_t|^2]\leq \cW_2^2(f_0,\tilde f_0)
\exp\left(C_\ga\intot (J_\ga(f_s)+J_\ga(\tilde f_s))ds\right),
$$
which concludes the proof since $\cW^2_2(f_t,\tilde f_t) \leq 
\E[|V_t-\tilde V_t|^2]$.
\end{preuve}

\section{Applications}\label{appli}
\setcounter{equation}{0}

We now want to show that the uniqueness result proved in the previous
section is relevant, in the sense that solutions in 
$L^\infty([0,T],\calp_2) \cap L^1([0,T],\cJ_\ga)$ indeed exist.
We will use Remark \ref{rk:Lp}, which says that 
$L^p\cap L^1 \subset \cJ_\ga $ for $p>3/(3+\ga)$.

We first recall that moments of solutions propagate, we give some
ellipticity estimate on the diffusion coefficient, and recall the chain
rule for the Landau equation.

Then Corollary \ref{thm:existence}-(i) is obtained by 
using the dissipation of entropy.
Finally, Corollary \ref{thm:existence}-(ii) 
is checked, using a direct computation.

\subsection{Moments}

We first recall the following result, which shows
that moments of the solutions to the Landau equation propagate.
We will use moments only in the case where 
$\gamma\in[-2,0)$, and the proposition below is quite easy.
We refer to \cite[Section 2.4. p 73]{villani:h}.
When $\gamma\in (-3,-2]$, the situation is much more delicate,
but Villani \cite[Appendix B p 193]{villani:these} has also proved
the propagation of moments.

\begin{prop}\label{prop:moments}
Let $\gamma \in(-2,0)$. Let us consider a weak solution 
$(f_t)_{t\in [0,T]}$ to (\ref{eq:Landau}).
Assume that for some $k\geq 2$, $m_{k}(f_0) \leq M<\infty$.
There is a constant $C$, depending only on $k,\ga,M,T$ 
such that 
$\sup_{[0,T]} m_{k}(f_s) \leq C$.
\end{prop}

\subsection{Ellipticity}
We need a result as Desvillettes-Villani
\cite[Proposition 4]{desvillettes-villani:00} 
(who work with $\gamma\geq 0$) on the ellipticity of the matrix $a$.

\begin{prop}\label{prop:ellipticity}
Let $\ga \in [-2,0)$. Let $E_0>0$ and $H_0>0$ be two constants, and
consider a nonnegative function $f$ such that 
${\int_{\R^3} f(v) dv=1}$, $m_2(f)\leq E_{0}$ and $H(f)\leq H_0$. 
There exists a constant $c>0$ depending on $\ga,E_{0}, H_0$ such that
\beas 
\forall \xi\in\R^3, \quad 
\sum_{i,j}\overline{a}^f_{ij}(v)\xi_i\xi_j\geq c(1+|v|)^\ga|\xi|^2
\eeas
where $\overline{a}^f(v)=\int_{\R^3} a(v-\vs)f(\vs)d\vs$.
\end{prop}

\begin{proof}
Since $\ga \in [-2,0]$, the proof of 
\cite[Proposition 4]{desvillettes-villani:00} can be applied: 
following line by line their proof, one can check that they use only that
$\gamma+2\geq 0$.
\end{proof}

One could easily extend this result to the case where $\ga \in (-3,-2)$.
However, we will not use it.

\subsection{Chain rule for the Landau equation}

As noted by Desvillettes-Villani \cite[Section 6]{desvillettes-villani:00},
we may write, for $f$ a weak solution to the Landau equation and
$\beta$ is a $C^1$ function with $\beta(0)=0$, at least formally,
\be\label{eq:beta}
\frac{d}{dt}\int_{\R^3}\beta(f_t(v))dv=-\int_{\R^3} \overline{a}(t,v)
\nabla f_t(v)\nabla f_t(v) \beta''(f_t)(v)dv-\int_{\R^3}\overline{c}(t,v)
\phi_\beta(f_t(v))dv
\ee
where 
\beas
& \overline{a}(t,v)=\overline{a}^{f_t}(v) =\int_{\R^3} a(|v-\vs|) 
f_t(\vs)d\vs, \\
& \overline{a}(t,v)\nabla f_t(v)\nabla f_t(v)=\sum_{i,j}
\overline{a}_{ij}(t,v)(\nabla f_t(v))_i(\nabla f_t(v))_j,\\
& \overline{c}(t,v)=-2(\ga+3)\int_{\R^3}|v-\vs|^\ga f_t(\vs)d\vs,\\
& \phi_\beta'(x)=x\beta''(x) \quad \text{and} \quad \phi_\beta(0)=0.
\eeas

\subsection{Moderately soft potentials}

Using the dissipation of entropy, we will deduce, at least for
$\ga$ not too much negative, the $L^p$ estimate we need.
Such an idea was handled in the much more delicate
case of the Boltzmann equation by
Alexandre-Desvillettes-Villani-Wennberg \cite{advw}.

\begin{prop}\label{prop:Lp-norm}
We assume that $\ga\in(-2,0)$. Let $\eps\in(0,1)$ with
$3-\eps>3/(3+\ga)$. Consider a weak solution 
$(f_t)_{t\in [0,T]}$ to (\ref{eq:Landau})
starting from $f_0$ with
$H(f_0)<\infty$, $m_2(f_0)<\infty$ and $m_q(f_0)<\infty$ with 
$q > 3|\ga|(2-\eps)/\eps$. Then, at least formally, 
$f \in L^1([0,T],L^{3-\eps})$.
\end{prop}

Before proving Proposition \ref{prop:Lp-norm}, we show how it allows
us to conclude the well-posedness for the Landau equation when
$\ga\in (-2,0)$.

\vip

\begin{preuve} {\it of Corollary \ref{thm:existence}-(i)}. 
We only have to prove the existence, since the uniqueness 
immediately from Theorem \ref{thm:uniqueness} and
Remark \ref{rk:Lp}.

We thus assume that $\ga \in (-2,0)$, and consider an initial condition
$f_0 \in \calp_2(\R^3)$, with $H(f_0)<\infty$, and
$m_2(f_0)<\infty$. Then Villani \cite{villani:98b} 
has shown the existence of a weak 
solution $(f_t)_{t\in [0,T]} \in L^\infty ([0,T],\calp_2)$, with
constant energy and nonnincreasing entropy, that is $m_2(f_t)=m_2(f_0)$
and $H(f_t) \leq H(f_0)$ for all $t\in [0,T]$.
We now use that $m_q(f_0)<\infty$, for some $q>q(\ga):=\ga^2/(2+\ga)$,
and we consider $p \in (3/(3+\ga),(3q-3\ga)/(q-3\ga))$.
Then one may check that for $\eps>0$ such that $p=3-\eps$, we have
$3|\ga|(2-\eps)/\eps < q$. Thus the {\it a priori} estimate 
proved in Proposition \ref{prop:Lp-norm} implies that 
the solution $(f_t)_{t\in [0,T]}$ can be built
in such a way that it lies in $L^1([0,T],L^p)$, which concludes the proof.
\end{preuve}

\vip

\begin{preuve} {\it of Proposition \ref{prop:Lp-norm}}.
We thus consider a weak solution $(f_t)_{t\in[0,T]}$ to the Landau equation.
Then this solution satisfies, at least formally, $H(f_t) \leq H(f_0)$
and $m_2(f_t)=m_2(f_0)$, for all $t\in [0,T]$. As a consequence, 
the ellipticity estimate of Proposition \ref{prop:ellipticity} 
is uniform in time when
applied to $f_t$.
We now divide the proof into several steps.

\vip

{\it Step 1.}
We apply (\ref{eq:beta}) with the function $\beta(x)=(x+1)\ln(x+1)$. 
One easily checks that $\beta''(x)=\frac{1}{x+1}$ and  
$0\leq \phi_\beta(x)=x-\ln(x+1)\leq x$.
Since $H(f_0)<\infty$ by assumption, we easily see that
${\int_{\R^3}\beta(f_0(v))dv<\infty}$. Using Proposition \ref{prop:ellipticity},
there exists a positive constant $c$ 
(depending only on $\ga,H(f_0),m_2(f_0)$) such that
\be \label{eq:global}
\frac{d}{dt}\int_{\R^3}\beta(f_t(v))dv\leq -c\int_{\R^3}(1+|v|)^\ga
\frac{|\nabla f_t(v)|^2}{1+f_t(v)}dv+2(\ga+3)\iint_{\R^3\times \R^3}
|v-\vs|^\ga f_t(v)f_t(v^*)dvdv^*.
\ee
First,
\beas
I_1&:=&\int_{\R^3}(1+|v|)^\ga\frac{|\nabla f_t(v)|^2}{1+f_t(v)}dv
=4 \int_{\R^3}(1+|v|)^\ga |\nabla(\sqrt{1+f_t(v)}-1)|^2 dv \\
&\geq& 4 (1+R)^\ga || \nabla(\sqrt{1+f_t}-1)||^2_{L^2(B_R)},
\eeas
for any $R>0$, where $B_R= \{x\in\R^3:|x|<R\}$.
By Sobolev embedding (see for example Adams \cite{adams:75}), we also know 
that there exists a constant $C>0$ such that for any $R>0$,
any measurable $g:\R^3\mapsto \R$ 
\beas
||g||_{L^6(B_R)}\leq C||g||_{H^1(B_R)}=C\left(||\nabla g||_{L^2(B_R)}+
||g||_{L^2(B_R)}\right).
\eeas
Consequently, there exists a constant $C>0$ such that
\begin{align*}
||\nabla(\sqrt{1+f_t}-1)&||^{2}_{L^2(B_R)}\geq C||\sqrt{1+f_t}-1||^2_{L^6(B_R)}-
||\sqrt{1+f_t}-1||^2_{L^2(B_R)}\\
&\geq C||\sqrt{f_t}\mathbf{1}_{\lbrace f_t\geq 1\rbrace}||^2_{L^6(B_R)}
-||f_t||_{L^1(B_R)}\geq C ||f_t\mathbf{1}_{\lbrace f_t\geq 1\rbrace}||_{L^3(B_R)}-
||f_t||_{L^1(B_R)}.
\end{align*}

Finally, since $||f_t||_{L^1(B_R)}\leq ||f_t||_{L^1(\R^3)}=1$, 
we get, for some $C>0$, for all $R\geq 1$ 
(which implies $(1+R)\leq 2R$), 
\begin{eqnarray}
I_1&\geq& 2^{2+\ga}R^{\ga}\left(C||f_t
\mathbf{1}_{\lbrace f_t\geq 1\rbrace}||_{L^3(B_R)}-1\right). \label{A}
\end{eqnarray}

Next we use Remark \ref{rk:Lp}, the H\"older inequality
and that $||f_t||_{L^1}=1$, to get,
for $p\in (3/(3+\ga),3-\eps)$, for some $C=C_{\ga,p}$,
\begin{eqnarray}
I_2 &:=& \iint_{\R^3\times\R^3} \!\!\!|v-\vs|^\ga f_t(v)f_t(v^*)dvdv^* 
\leq ||f_t||_{L^1} J_\ga(f_t) \nonumber \\
&\leq& C
\left(1+||f_t||_{L^p}\right) 
= C\left(1+\left(\int_{\R^3} f_t(v) f_t^{p-1}(v) dv
\right)^{\frac{1}{p}}\right) \nonumber \\
&\leq& C\left(1+\left(\int_{\R^3}f_t(v) f_t^{2-\eps}(v)dv
\right)^{\frac{(p-1)}{(2-\eps)p}}\right)=C\left(1+||f_t||_{L^{3-\eps}}^{\frac{(p-1)
(3-\eps)}{(2-\eps)p}}
\right). \label{B}
\end{eqnarray}
Set $\delta=\frac{(p-1)(3-\eps)}{(2-\eps)p}$. Using 
(\ref{eq:global}-\ref{A}-\ref{B}), we deduce that there is $C>0$ such that for 
any $R\geq1$,
\begin{eqnarray*}
 R^{\ga}\int_0^T||f_s\mathbf{1}_{\lbrace f_s\geq 1\rbrace}||_{L^3(B_R)}ds &\leq& 
C\left(\int_{\R^3}\beta(f_0(v))dv-\int_{\R^3}\beta(f_T(v))
dv+\int_0^T \left(1+|| f_s||_{L^{3-\eps}}^{\delta}\right)ds\right) \\
&\leq& C\left(1+ \int_0^T ||f_s||_{L^{3-\eps}}^{\delta} 
ds\right).
\end{eqnarray*}

{\it Step 2.} For $\alpha>0$, we define  $g_t(v)=f_t(v)(1+|v|)^{\ga-\al}
\mathbf{1}_{\lbrace f_t\geq 1\rbrace}$.
We have, using Step 1,
\beas
\int_0^T||g_s||_{L^3}ds &\leq& 
\int_0^T\sum_{k\geq 0}
||g_s||_{L^3(\{2^k -1\leq|v|\leq2^{k+1} -1\})}ds 
\leq \int_0^T \sum_{k\geq 0}2^{k(\ga-\al)}
||f_s\mathbf{1}_{\lbrace f_s\geq 1\rbrace}||_{L^3(B_{2^{k+1}})}ds
\nonumber\\
&\leq& C2^{-\ga}\sum_{k\geq 0}2^{-\al k}
\left( 1+ \int_0^T || f_s||_{L^{3-\eps}}^{\delta} ds\right) 
\leq C\left( 1+ \int_0^T || f_s||_{L^{3-\eps}}^{\delta} ds\right).
\eeas

\vip

{\it Step 3.} We now prove if $\alpha>0$ is small enough,
for some constant $C$,
\be\label{eq:g}
||f_t\mathbf{1}_{\lbrace f_t\geq 1\rbrace}||_{L^{3-\eps}}\leq C\left(1+ 
||g_t||_{L^3}\right).
\ee
We consider a nonnegative function $h$ with $\int_{\R^3}h(v)dv
\leq 1$. By H\"older's inequality, for 
$\eps\in(0,1)$
\beas
||h||_{L^{3-\eps}}&=&\left(\int_{\R^3}h^{2-\eps}(v)
\left(\frac{(1+|v|)^{(\ga-\al)}}{(1+|v|)^{(\ga-\al)}}\right)^{3(2-\eps)/2}h(v)dv
\right)^{1/(3-\eps)}\\
&\leq&\left(\int_{\R^3}\left[(1+|v|)^{\ga-\al}h(v)\right]^{3}dv
\right)^{\frac{2-\eps}{2(3-\eps)}}\left(\int_{\R^3}(1+|v|)^{{3(\al-\ga)
(2-\eps)}/{\eps}}h(v)dv\right)^{\frac{\eps}{2(3-\eps)}}\\
&\leq&\left(\int_{\R^3}\left[(1+|v|)^{\ga-\al}h(v)\right]^{3}dv\right)^{
\frac{2-\eps}{2(3-\eps)}}\left(1+\int_{\R^3}|v|^{{3(\al-\ga)(2-\eps)}/{\eps}}
h(v)dv\right).
\eeas
Then, for $h=f_t\mathbf{1}_{\lbrace f_t\geq 1\rbrace}$, setting 
$r=3(\al-\ga)(2-\eps)/\eps$, and recalling that 
$g_t(v)=f_t(v)(1+|v|)^{\ga-\al}\mathbf{1}_{\lbrace f_t\geq 1\rbrace}$
\beas
||f_t\mathbf{1}_{\lbrace f_t\geq 1\rbrace}||_{L^{3-\eps}}\leq\left(1+m_{r}(f_t)
\right) ||g_t||_{L^3}^{\frac{3(2-\eps)}{2(3-\eps)}}\leq(1+m_{r}(f_t)) (1+
||g_t||_{L^3}).
\eeas
But by assumption, $m_q(f_0)<\infty$ for some $q>3|\ga|(2-\eps)/\eps$,
whence, by Proposition \ref{prop:moments}, 
$\sup_{[0,T]} m_q(f_t) <\infty$. Choosing $\alpha>0$ small
enough, in order that $q\geq r$, we deduce (\ref{eq:g}).

\vip

{\it Step 4.} Using that $\int_{\R^3} f_s(v)dv=1$, 
Steps 2 and 3, we obtain, for some constant $C$
(depending in particular on $T$),
$$
\int_0^T ||f_s||_{L^{3-\eps}}ds \leq C+ 
\int_0^T ||f_s\mathbf{1}_{\lbrace f_s\geq 1\rbrace}||_{L^{3-\eps}}ds
\leq C +C \int_0^T ||f_s||_{L^{3-\eps}}^\delta ds
\leq C + C \left( \int_0^T ||f_s||_{L^{3-\eps}} \right)^\delta .
$$
But one may choose $p\in (3/(3+\ga),3-\eps)$ (recall Step 1) such that 
$\delta=\frac{(p-1)(3-\eps)}{p(2-\eps)} \in (0,1)$ (choose $p$
very close to $3/(3+\ga)$ and use that by assumption,
$3/(3+\ga)< 3-\eps$, whence $\eps<\frac{6+3\ga}{3+\ga}$).
As a consequence, $\int_0^T ||f_s||_{L^{3-\eps}}ds \leq x_0$,
where $x_0$ is the largest solution of $x=C+Cx^\delta$.
Following carefully the proof above, one may check that
$C$, and thus $x_0$, depends only on $T,f_0,q,\gamma,\eps$.
\end{preuve}

\subsection{Soft potentials}

We now would like to obtain a result which includes the case
of very soft potentials, that is $\ga \in(-3,-2]$.

\begin{prop}\label{prop:Lp-norm2}
Let $\ga\in(-3,0)$ and $p\in(\frac{3}{3+\ga},+\infty)$. 
Let $f$ be a weak solution of the Landau equation
with $f_0\in \calp_2(\R^3)\cap L^p(\R^3)$. 
Then there exists a time $T^*>0$ depending on $\ga$, $p$ 
and $||f_0||_{L^p}$ such that for any $T\in[0,T^*)$, at
least formally,
$\sup_{[0,T]}||f_t||_{L^p}<\infty$.
\end{prop}

\begin{proof}
Let us consider the function $\beta(x)=x^p$. Since 
$p>1$, we have $\beta''\geq 0$ and 
$\phi_\beta(x)=(p-1)x^p$. Using (\ref{eq:beta}), neglecting all
nonnegative terms, and using Remark 
\ref{rk:Lp}, since $p>\frac{3}{3+\ga}$, there exists a constant 
$C_{\ga,p}>0$ such that
\beas
\frac{d}{dt}||f_t||_{L^p}^p&\leq& 2(\ga+3)(p-1)
\iint_{\R^3\times \R^3 }|v-\vs|^\ga f_t(\vs)f^p_t(v)dvd\vs\\
&\leq&||f_t||_{L^p}^p J_\ga(f_t) \leq C_{\ga,p}(1+||f_t||_{L^p})||f_t||_{L^p}^p
\leq C_{\ga,p}(1+||f_t||_{L^p}^{2p}) .
\eeas
Thus for all 
$0\leq t<T^*:=(\frac{\pi}{2}-\arctan||f_0||_{L^p}^p)/C_{\ga,p}$, we have
$||f_t||^p_{L^p}\leq \tan(\arctan||f_0||^p_{L^p}+C_{\ga,p}t)$,
which concludes the proof.
\end{proof}

\vip

\begin{preuve} {\it of Corollary \ref{thm:existence}-(ii)}. 
We only have to prove the existence, since the uniqueness follows 
immediately follows from Theorem \ref{thm:uniqueness} and
Remark \ref{rk:Lp}.

We thus assume that $\ga \in (-3,0)$, $p>3/(3+\ga)$, 
and consider an initial condition
$f_0 \in \calp_2(\R^3) \cap L^p(\R^3)$ (which implies that 
$H(f_0)<\infty$). Then Villani \cite{villani:98b} 
has shown the existence of a weak 
solution $(f_t)_{t\in [0,T]} \in L^\infty ([0,T],\calp_2)$,
for arbitrary $T$.

But using the {\it a priori} estimate of Proposition \ref{prop:Lp-norm2}, 
we deduce that this solution can be built in such a way that it
belongs to $L^\infty_{loc}([0,T_*),L^p)$. This concludes the proof.
\end{preuve}

\end{document}